\documentclass[12pt]{article}

\usepackage{pst-poly}     
\usepackage{pst-all}      
\usepackage{amsfonts,amssymb,amsmath,latexsym,bm,array}
\usepackage{cite}
\newtheorem{thm}{Theorem}[section]
\newtheorem{Def}[thm]{Definition}

\newtheorem{cor}[thm]{Corollary}
\newtheorem{lem}[thm]{Lemma}
\newtheorem{cjt}[thm]{Conjecture}
\newtheorem{pro}[thm]{Problem}
\newenvironment{pf}[1][Proof]{\noindent\textbf{#1.} }{\hfill\rule{1mm}{2mm}}

\makeatletter \@addtoreset{equation}{section} \makeatother

\begin{document}
\title{Generalization of the cover pebbling number on trees\footnote{The research of Zheng-Jiang Xia is supported by Key Projects in Natural Science Research of Anhui Provincial Department of Education (No. KJ2018A0438). The research of Zhen-Mu Hong is supported by NSFC (No.11601002).}
}
\author
{Zheng-Jiang Xia\thanks{Corresponding author, Email: xzj@mail.ustc.edu.cn (Z.-J. Xia), zmhong@mail.ustc.edu.cn
(Z.-M. Hong).},\quad Zhen-Mu Hong
\\ \\
{\small     School of Finance,} \\
{\small           Anhui University of Finance \& Economics,}   \\
{\small             Bengbu, Anhui, 230030, P. R. China}\\
}

\date{}
\maketitle

\begin{quotation}
\noindent\textbf{Abstract}: A pebbling move on a graph consists of taking two pebbles off from one vertex and add one pebble on an adjacent vertex, the $t$-pebbling number of a graph $G$ is the minimum number of pebbles so that we can move $t$ pebbles on any vertex on $G$ regardless the original distribution of pebbles. Let $\omega$ be a positive function on $V(G)$, the $\omega$-cover pebbling number of a graph $G$ is the minimum number of pebbles so that we can reach a distribution with at least $\omega(v)$ pebbles on $v$ for all $v\in V(G)$. In this paper, we give the $\omega$-cover pebbling number of trees for nonnegative function $\omega$, which generalized the $t$-pebbling number and the traditional weighted cover pebbling number of trees.



\end{quotation}

\section{Introduction}

Pebbling in graphs was first introduced by Chung\cite{c89}. For a given connected graph $G=(V,E)$, a \emph{distribution} $D$ of $G$ is a projection from $V(G)$ to the nonnegative integers, $D(v)$ represents the number of pebbles on the vertex $v$, the total number of pebbles on a subset $A$ of $V$ is given by $|D(A)|=\sum_{v\in A}D(v)$, $|D|=|D(V)|$ is the size of $D$. $d(u,v)$ is the distance of $u$ and $v$, and we write $u\sim v$  if they are adjacent. $N(v)=\{u|u\sim v\}$ is the neighbor of $v$, $d(v)=|N(v)|$ is the degree of $v$, let $H$ be an induced subgraph of $G$, we use $d_H(v)$ to denote the degree of $v$ in $H$.

\begin{Def}\rm
A \emph{pebbling move} consists of the removal of two pebbles from a vertex and the placement of one pebble on an adjacent vertex.
\end{Def}

\begin{Def}\rm
\emph{The $t$-pebbling number} of a vertex $v$ in $G$, denoted by $f_t(G,v)$, is the minimum number of pebbles that are sufficient to move $t$ pebbles to $v$ regardless of the original distribution of pebbles. \emph{The $t$-pebbling number }of $G$, $f_t(G)=\max_{v\in V(G)}f_t(G,v)$. \emph{The pebbling number} of $G$ is $f_1(G)$, and we denote it $f(G)$.
\end{Def}

To determine the pebbling number of general graph is difficult. The problem of whether a distribution can reach a fixed vertex was shown to be NP-complete\cite{lcd14,mc06}. The problem of deciding if the pebbling number of a graph $G$ is less than $k$ was shown to be $\Pi_2^P$-complete\cite{mc06}. The pebbling numbers of trees\cite{m92}, cycles\cite{psv95}, hypercubes\cite{c89}, squares of cycles\cite{yzz12,yz12} and so on have been determined. There is a conjecture given by Chung\cite{c89}, which is called Graham's Conjecture.

\begin{cjt}{\rm(Graham's Conjecture)}
Let $G$ and $H$ be two connected graphs, the pebbling number of the Cartesian product of $G$ and $H$ satisfies:
$$f(G\times H)\leq f(G)f(H).$$
\end{cjt}

There are many results about Graham's Conjecture\cite{h05}, but this conjecture is still open.

We first introduce path partition and the pebbling number of trees.

\begin{Def}\rm{(\cite{m92})}
Given a root vertex $v$ of a tree $T$, then we can view $T$ be a directed graph $\overrightarrow{T_v}$ with each edge directed to $v$, \emph{a path
partition} is a set of nonoverlapping directed paths the union of which is  $\overrightarrow{T_v}$. A path partition
is said to \emph{majorize }another if the nonincreasing sequence of the path size majorizes
that of the other (that is $(a_1,a_2,\ldots,a_r)>(b_1,b_2,\ldots,b_t)$ if and only if $a_i > b_i$  where $i= \min\{j:a_j\neq b_j \}$). A path partition of a tree $\overrightarrow{T_v}$ is said to be \emph{maximum} if it majorizes all other path partitions.
\end{Def}

\textbf{Note:} By the definition of the maximum path partition, we can give a way to determine the size of the maximum path partition: first we choose a longest directed path $P_1$ in $\overrightarrow{T_v}$, with length $a_1$, then we choose a longest directed path $P_2$ in $\overrightarrow{T_v}\backslash E(P_1)$, with length $a_2$, and so on.

Moews\cite{m92} found the $t$-pebbling number of trees by a path partition.
\begin{thm}{\rm(\cite{m92})}
Let $T$ be a tree, $v\in V(T)$, $(a_1,\ldots, a_n)$ is the size of the maximum path partition of $\overrightarrow{T_v}$, then
$$f_t(T,v)=t2^{a_1}+\sum_{i=2}^{n}2^{a_i}-n+1,$$
$$f_t(T)=\max_{v\in V(T)}f_t(T,v).$$
\label{t1}\end{thm}

\begin{cor}
Let $T$ be a tree, $v\in V(T)$, $\alpha=(a_1,\ldots, a_n)$ is the size of a path partition of $\overrightarrow{T_v}$, $s_\alpha:=t2^{a_1}+\sum_{i=2}^{n}2^{a_i}-n+1$, then
$$f_t(T,v)=\max_\alpha s_\alpha.$$
\label{cor0}\end{cor}

\begin{pf}
Let $\alpha_0$ be the size of the maximum path partition of $\overrightarrow{T_v}$, then $f_t(T,v)=s_{\alpha_0}\leq \max_\alpha s_\alpha.$

Let $P_1,P_2,\ldots,P_n$ be a path partition of $\overrightarrow{T_v}$, and the length of $P_i$ is $a_i$ for $1\leq i\leq n$. Note that for each $P_i$, assume the two endpoints $v_i$ and $v_i'$ satisfies $d(v_i,v)>d(v_i',v)$. We put $t2^{a_1}-1$ pebbles on $v_1$, and $2^{a_i}-1$ pebbles on $v_i$ for $2\leq i\leq n$, it is clear that we cannot move $t$ pebbles on $v$ from this distribution. Thus for each $\alpha$, $s_\alpha-1<f_t(T,v)$, so $s_\alpha\leq f_t(T,v)$, so $\max_\alpha s_\alpha\leq f_t(T,v)$.
\end{pf}

\begin{Def}\rm
Let $\omega$ be a nonnegative function on $V(G)$, $D$ is a distribution on $V(G)$, we say $D$ is \emph{$\omega$-solvable} (or $D$ solves $\omega$), if we can reach a distribution $D^*$ from $D$, by a sequence of pebbling moves, so that $D^*(v)\geq \omega(v)$. The\emph{ $\omega$-cover pebbling number} of $G$, denoted by $\gamma_\omega(G)$, is the minimum number $\gamma_\omega(G)$ so that every distribution $D$ with $\gamma_\omega(G)$ pebbles is $\omega$-solvable.
\end{Def}

\begin{Def}\rm
Let $\omega$ be a positive function on $V(G)$, define
$$s_\omega(v)=\sum_{u\in V(G)}\omega(u)2^{d(u,v)},$$
and
$$s_\omega(G)=\max_{v\in V(G)}s_\omega(v).$$
\label{def1}\end{Def}

The $\omega$-cover pebbling number of a graph $G$ has been determined for \emph{positive} $\omega$ by \cite{s05,vw04}.

\begin{thm}{\rm(\cite{s05,vw04})}
Let $\omega$ be a positive weight function on $V(G)$, the $\omega$-cover pebbling number of $G$ is
$$\gamma_{\omega}(G)=s_\omega(G).$$
\label{t2}\end{thm}

From Theorem~\ref{t2}, one can show

\begin{thm}({\rm\cite{s05,vw04}})
Let $\omega_1$ be a positive function on $G$, and $\omega_2$ be a positive function on $H$, the function $\omega$ on $G\times H$ is given by $\omega((g,h))=\omega_1(g)\omega_2(h)$, where $g\in V(G)$ and $h\in V(H)$, then $\gamma_\omega(G\times H)= \gamma_{\omega_1}(G)\gamma_{\omega_2}(H)$.
\end{thm}

We first generalized the definition of $s_\omega(T)$ while $\omega$ is a \emph{nonnegative} function on a tree $T$.

\begin{Def}\rm
Given a tree $T$ and a nonnegative function $\omega$, for each vertex $v\in V(T)$, let $T_\omega(v)$ be the minimum subtree containing $v$ and $W:=\{u: \omega(u)>0\}$, and we give each edge in $T\backslash E(T_\omega(v))$ a direction towards $T_\omega(v)$, which is denoted by $\overrightarrow{T}\backslash E(T_\omega(v))$, and $(a_1,\ldots, a_n)$ is the size of the maximum path partition of $\overrightarrow{T}\backslash E(T_\omega(v))$. We define
$$s_\omega(v)=\sum_{u\in W}\omega(u)2^{d(u,v)}+\sum_{i=1}^{n}2^{a_i}-n.$$
and
$$s_\omega(T)=\max_{v\in V(T)}s_\omega(v).$$
\label{def2}\end{Def}

Note that while $\omega$ is positive, then the two definitions of $s_\omega(T)$ are the same, so Definition~\ref{def2} is a generalization of Definition~\ref{def1}.
We generalized Theorem~\ref{t2} while $T$ is a tree and $\omega$ is nonnegative, our main result in this paper is as follows:
\begin{thm}
Let $T$ be a tree with a nonnegative weight function $\omega$ on $V(T)$, the $\omega$-cover pebbling number of $T$ is
$$\gamma_{\omega}(T)=s_\omega(T).$$
\label{thmmain}\end{thm}

\begin{thm}
Let $T$ be a tree with a nonnegative weight function $\omega$ on $V(T)$, if $|W|=1$, then $Theorem~\ref{thmmain}\Leftrightarrow Theorem~\ref{t1}$.
\label{thm0}\end{thm}

\begin{pf}
If $|W|=1$, assume that $\omega(v)=t$, and $\omega(u)=0$ for $u\neq v$, we will show that $f_t(T,v)=s_\omega(T)$.

Assume the size of a maximum path partition of $\vec{T}_v$ is $(a_0,a_1,\ldots,a_n)$, and $d(v,v_0)=a_0$, $P_0$ be the longest directed path from $v_0$ to $v$. Then $(a_1,\ldots,a_n)$ must be the size of a maximum path partition in $\vec{T}_v\backslash P_0$. So $f_t(T,v)=s_\omega(v_0)\leq s_\omega(T)$.

Assume $s_\omega(T)=s_\omega(v_1)$, and $d(v_1,v)=a_0$, Let $P_0$ be the path connected $v_1$ and $v$, then $T_\omega(v_1)=P_0$, assume $(a_1,\ldots, a_n)$ is the size of the maximum path partition of $T\backslash E(T_\omega(v))=T\backslash E(P_0)$, so $\alpha=(a_0,a_1,\ldots, a_n)$ is a path partition of $\vec{T}_v$, and $s_\alpha=s_\omega(v_1)$, by Corollary~\ref{cor0}, $f_t(T,v)\geq s_\omega(v_1)=s_\omega(T)$, and it's over.
\end{pf}

\begin{Def}\rm(\cite{ccfh05})
Given a sequence $S$ of pebbling moves on $G$, \emph{the transition digraph} obtained from $S$ is a directed multigraph denoted $T(G,S)$ that has $V(G)$ as its vertex set, and each move $s\in S$ along edge $uv$ (move off two pebbles from $u$ and add one on $v$) is represented by a directed edge $uv$.
\end{Def}

The following lemma is useful in next sections.
\begin{lem}{\rm(\cite{ccfh05}, No-Cycle Lemma)}
Let $S$ be a sequence of pebbling moves on $G$, reaching a distribution $D$. Then there exists a sequence $S^*$ of pebbling moves, reaching a distribution $D^*$, such that

1. On each vertex $v$, $D^*(v)\geq D(v)$;

2. $T(G,S^*)$ does not contain any directed cycles.
\label{ncl}\end{lem}

\section{Preliminaries}

%
\begin{Def}\rm
Let $C$ be a \emph{generalized distribution} on $G$, satisfies $C(v)$ is an integer (may be negative) for all $v\in V(G)$. \emph{A pebbling move} on $G$  consists of the removal of two pebbles from a vertex $v$ (with $C(v)\geq2$) and the placement of one pebble on an adjacent vertex.
\end{Def}

In the following of this paper, a \emph{distribution} $D$ means that $D(v)\geq0$, and a \emph{generalized distribution} $C$ means $C(v)$ is an integer for all $v\in V(G)$.

\begin{Def}\rm
Let $\omega$ be a nonnegative function on $V(G)$, $C$ is a generalized distribution on $V(G)$, we say $C$ is \emph{$\omega$-solvable}, if we can reach a distribution $C^*$ from $C$, by a sequence of pebbling moves, so that $C^*(v)\geq \omega(v)$. In particular, if $\omega(v)=0$ for all $v\in V(G)$, then we say that $C$ is \emph{$0$-solvable}.
\end{Def}

\begin{lem}
Let $D$ be a distribution on a graph $G$ and $\omega$ be a nonnegative function on $V(G)$, $C:=D-\omega$. Then $D$ is $\omega$-solvable. $\Leftrightarrow$ $C$ is $0$-solvable.
\label{lem1}\end{lem}

\begin{pf}
If $C$ is $0$-solvable, then with the same sequence of pebbling moves, we can find that $D$ is $\omega$-solvable.

On the other side, if $D$ is $\omega$-solvable, by Lemma~\ref{ncl}, there exist a sequence of pebbling moves $S$ reaching a distribution $D^*$ with $D^*(v)\geq \omega(v)$ and $T(G,S)$ does not contain any direct cycle. So we can give a sequence of the vertices of $G$, as $(v_1,v_2,\ldots,v_n)$, so that each directed edge $v_iv_j$ in $T(G,S)$ satisfies $i<j$. Thus we can rearrangement the sequence of pebbling moves $S$ along the order $(v_1,v_2,\ldots,v_n)$, that means we first choose all pebbling moves in $S$ that move pebbles off $v_1$, then choose all pebbling moves in $S$ that move pebbles off $v_2$ and so on, denote this sequence of pebbling moves $S'$. Since no directed edges is from $v_j$ to $v_i$ for $j>i$, while we begin to move pebbles off $v_{i+1}$, then the number of pebbles left on $v_i$ is just $D^*(v_i)(\geq \omega(v_i)$), that means $S'$ is a sequence of pebbling moves reaching $D^*-\omega\geq0$ from $C$, thus $C$ is $0$-solvable.
\end{pf}

\begin{Def}\rm
Let $D$ be a distribution  on a tree $T$, $\omega$ is a nonnegative function on $V(T)$, $C=D-\omega$ is called the \emph{induced generalized distribution} from $D$ and $\omega$ of $T$. Let $v$ be a leaf of $T$, $u$ is its neighbor in $T$, the \emph{induced generalized distribution} $C'$ on $T\setminus v$ is given as follows: if $C(v)\geq0$, then $C'(u)=C(u)+\left\lfloor C(v)/2\right\rfloor$, if $C(v)<0$, then $C'(u)=C(u)+2C(v)$, keeping $C'(x)=C(x)$ unchanged for all $x\neq u$.
\end{Def}

\begin{lem}
Let $D$ be a distribution  on a tree $T$, $\omega$ is a nonnegative function on $V(T)$, $C:=D-\omega$, $v$ is a leaf of $T$, $C'$ is the induced generalized distribution from $D$ and $\omega$ of $T\backslash v$. Then $C$ is $0$-solvable in $T$. $\Leftrightarrow$ $C'$ is $0$-solvable in $T\backslash v$.
\label{lem2}\end{lem}

\begin{pf}
Firstly, we assume $C$ is $0$-solvable in $T$, and there is a sequence of pebbling moves $\sigma$ reaching a distribution $C^*$ from $C$ with $C^*(x)\geq0$ for each $x\in V(T)$.

Case $1.1$. $C(v)\geq 0$. By Lemma~\ref{ncl}, we may assume that no pebble has been moved from $u$ to $v$, so at most $\left\lfloor C(v)/2\right\rfloor$ pebbles can be moved from $v$ to $u$. We may assume the first step of $\sigma$ is to move $\left\lfloor C(v)/2\right\rfloor$ pebbles from $v$ to $u$, so the left steps makes $C'$ solves $0$ on $T\backslash v$, and we are done.

Case $1.2$. $C(v)<0$. By Lemma~\ref{ncl}, we may assume that no pebble has been moved from $v$ to $u$. So we may assume the last step of $\sigma$ is to move $-C(v)$ pebbles from $u$ to $v$, so the steps before it makes $C'$ solves $0$ on $T\backslash v$, and we are done.

Secondly, we assume $C'$ is $0$-solvable in $T\backslash v$, and there is a sequence of pebbling moves $\delta$ reaching a distribution $C^*$ from $C'$ with $C^*(x)\geq0$ for each $x\in V(T\backslash v)$.

Case $2.1$. $C(v)\geq 0$. First, we move $\left\lfloor C(v)/2\right\rfloor$ pebbles from $v$ to $u$, and the left steps are just $\delta$, this sequence makes $C$ solves $0$.

Case $2.2$. $C(v)<0$. After the sequence of pebbling moves $\delta$, we move $-C(v)$ pebbles from $u$ to $v$, this sequence makes $C$ solves $0$, over.
\end{pf}

\textbf{Notations:} Assume $T^*$ is a subtree of $T$, then $T^*$ can be obtained from $T$ by deleting the leaf of the subtree of $T$ (the vertex with degree one) one by one, so for each subtree $T^*$ of $T$, we can get an induced generalized distribution $C^*$. In particular, for each vertex $v\in V(T)$, let $T_v$ be a subtree containing $v$ and all of its neighbors. We use $C_v$ to denote the induced generalized distribution from $D$ and $\omega$ of $T_v$, and $\widehat{C}(v)$ to denote the induced generalized distribution of $\{v\}$.

\begin{cor}
Let $D$ be a distribution  on a tree $T$, $\omega$ is a nonnegative function on $V(T)$, and $\widehat{C}(v)$ is the induced generalized distribution from $D$ and $\omega$ of $\{v\}$. $D$ is not $\omega$-solvable.$\Leftrightarrow$ $\widehat{C}(v)<0$ for each $v\in V(T)$.
\label{cor3}\end{cor}

\begin{pf}
From Lemma~\ref{lem1} and Lemma~\ref{lem2}, the result follows by induction.
\end{pf}

\begin{lem}
Let $D$ be a distribution on a tree $T$ which is not $\omega$-solvable with $|D|=\gamma_\omega(T)-1$, then for each vertex $x\in V(T)$ which is not a leaf of $T$, there exist a vertex $y\in N(x)$, so that $C_x(y)\geq0$.
\label{lem3}\end{lem}

\begin{pf}
If $C_x(x')< 0$, for all $x'\in N(x)$. Assume $y,z\in N(x)$ satisfies $C_x(z)\leq C_x(y)<0$, then we delete all other vertices to left $T_1=yxz$, and its induced generalized distribution $C_1$. Then $C_1(y)=C_x(y)$, $C_1(z)=C_x(z)$, and $\widehat{C}(x)=C_1(x)+2C_1(y)+2C_1(z)\leq -1$ by Corollary~\ref{cor3}. Note that $C_1(x)=D(x)-w(x)+\sum_{x'\in N(x),x'\notin \{y,z\}}2C_x(x')$. So $C_1(x)-D(x)\leq 0$ and $C_1(x)+2C_1(z)-D(x)\leq 0$.
Now we remove $D(x)$ pebbles from $x$, and put $D(x)+1$ pebbles on $y$ to get a new distribution $D'$ with $|D'|=|D|+1$, the induced generalized distribution from $D'$ and $\omega$ of $\{y\}$ is denoted by $\widehat{C'}(y)$. Then $\widehat{C'}(y)=(C_1(y)+D(x)+1)+2(C_1(x)+2C_1(z)-D(x))=(C_1(x)+2C_1(y)+2C_1(z))+(C_1(z)-C_1(y))+C_1(z)+(C_1(x)-D(x))+1\leq -1+0-1+0+1=-1$, so $D'$ is not $\omega$-solvable by Corollary~\ref{cor3}, a contradiction to $|D'|=\gamma_\omega(T)$, and we are done.
\end{pf}

\begin{lem}
Let $D$ be a distribution which is not $\omega$-solvable in $T$, $x\sim y$. if $C_x(y)\geq0$, then $C_y(x)<0$.
\label{lem4}\end{lem}

\begin{pf}
Let $T_1=xy$, and its induced generalized distribution is denoted by $C_1$, then $C_1(x)=C_y(x)$, $C_1(y)=C_x(y)$, if both of them are nonnegative, then $\widehat{C}(x)=C_1(x)+\left\lfloor C_1(y)/2\right\rfloor\geq 0$. By Corollary~\ref{cor3}, $D$ is $\omega$-solvable, a contradiction to the condition $D$ is not $\omega$-solvable, and we are done.
\end{pf}

\begin{thm}
Let $\omega$ be a nonnegative function on $V(T)$,  there exist a distribution $D$, which is not $\omega$-solvable with $|D|=\gamma_\omega(T)-1$, and all pebbles are distributed on the leaves of $T$.
\label{thm1}\end{thm}

\begin{pf}
We will construct such distribution as follows. Assume that $D$ is not $\omega$-solvable with $|D|=\gamma_\omega(T)-1$, $x$ is a vertex which is not a leaf, and $D(x)>0$, then $\widehat{C}(x)<0$ from Corollary~\ref{cor3}.

From Lemma~\ref{lem3}, we may assume that there exist a vertex $y\in N(x)$, with $C_x(y)\geq 0$, the induced generalized distribution of $T_1:=xy$ is denoted by $C_1$. Then $C_1(y)=C_x(y)\geq0$, $\widehat{C}(x)=C_1(x)+\left\lfloor C_1(y)/2\right\rfloor<0$. Now we move all pebbles from $x$ to $y$ to get a new distribution $D'$ with $|D'|=|D|$, the induced generalized distribution from $D'$ and $\omega$ of $y$ is denoted by $\widehat{C'}(y)$. $\widehat{C'}(y)=(C_1(y)+D(x))+2(C_1(x)-D(x))<0$, so $D'$ is not $\omega$-solvable by Corollary~\ref{cor3}.

Now we consider the new distribution $D'$ on $T$. We use $C'_v$ to denote the induced generalized distribution from $D'$ and $\omega$ of $T_v$ for $v\in V(T)$. From Lemma~\ref{lem4}, $C'_y(x)<0$, by Lemma~\ref{lem3}, there must exist $z\in N(y)$ so that $C'_y(z)\geq0$, so we can remove $D'(y)$ pebbles from $y$ and put $D'(y)$ pebbles on $z$ to get a new distribution $D''$ which is not $\omega$-solvable and so on, until we move the pebbles from $x$ to some leaf of $T$, and we can do the same thing to other vertex $x'$ with $D(x')>0$, and we are done.
\end{pf}

\section{The generalization of the cover pebbling number on trees}

Assume that $s_\omega(v_0)=s_\omega(T)$ for some $v_0\in V(T)$, note that $\overrightarrow{T}\backslash E(T_\omega(v_0))$ is a directed graph, we define $d_\omega(u,l)$ be the length of the maximal path containing $u$ in all maximum path partitions of $\overrightarrow{T}\backslash E(T_\omega(v_0))$, if $\omega$ is clear, then we use $d(u,l)$ for short (note that $d(u,l)$ maybe $0$). Let $P_\alpha$ be a maximal path partition of  $\overrightarrow{T}\backslash E(T_\omega(v_0))$, Then $d_\omega(u,l)=\max_{P_\alpha}\{|P|:u\in P,P\in P_\alpha\}$.

%

\begin{lem}
Assume that $s_\omega(v_0)=s_\omega(T)$ for some $v_0\in V(T)$, then for each vertex $u\in V(T)$, $d(u,v_0)\geq d(u,l)$.
\label{lem5}\end{lem}

\begin{pf}
If $|W|=1$, we may assume that $\omega(v)=t$, and $\omega(u)=0$ for $u\neq v$. By the proof of Theorem~\ref{thm0}, we know that $f_t(T,v)=s_\omega(v_0)$. Let $(a_1,a_2,\ldots,a_n)$ be the size of the maximum path partition of $\overrightarrow{T_v}$. Then $d(v,v_0)=\max_{u\in V(T)} d(v,u)=a_1$. Assume $P_1$ be the path connected $v$ and $v_0$. $P_2$ be the maximal path containing $u$ in $\overrightarrow{T_v}\backslash P_1$, and $P_1\cap P_2=v'$, and the path connected $v$ and $v'$ is denoted by $P_3$. Then the length of $P_1$ ($P_2$) is $a_1$ ($d(u,l)$), and $d(v',v_0)\leq d(u,v_0)$. If  $d(u,v_0)< d(u,l)$, then $d(v',v_0)< d(u,l)$, we get a path $P_2\cup P_3$ with length $a_1-d(v',v_0)+d(u,l)>a_1$, a contradiction to the maximum of $a_1$, thus $d(u,v_0)\geq d(u,l)$.

If $|W|\geq 2$. We only need to show it while $u\in V(T_\omega(v_0))$.

If $d(u,v_0) < d(u,l)$ for some $u\in V(T_\omega(v_0))$, let $u'$ be some vertex in $W\backslash u$.
Then there exist some leaf $v_1$ in $\overrightarrow{T}\backslash E(T_\omega(v_0))$ so that $d(u,v_0)<d(u,v_1)$, we will show that $s_\omega(v_1)>s_\omega(v_0)$.

Let $P$ be the path connected $v_1$ and $u$, then we know the inner vertices of $P$ do not belong to $W$, so for each $x\in W\backslash u$, $d(x,v_1)>d(u,v_1)$ and $d(x,v_1)-d(x,v_0)\geq d(u,v_1)-d(u,v_0)+1\geq2$.

Note that $\overrightarrow{T}\backslash E(T_\omega(v_0))\backslash P\subseteq \overrightarrow{T}\backslash E(T_\omega(v_1))$. So
\begin{align*}
&s_\omega(v_1)-s_\omega(v_0)\\
&\geq\sum_{x\in W}\omega(x)(2^{d(x,v_1)}-2^{d(x,v_0)})-2^{d(u,v_1)}\\
&\geq \omega(u')(2^{d(u',v_1)}-2^{d(u',v_0)})-2^{d(u,v_1)}\\
&\geq2^{d(u',v_1)}-2^{d(u',v_0)}-2^{d(u,v_1)}\\
&\geq2^{d(u',v_1)}-\frac{2^{d(u',v_1)}}{4}-\frac{2^{d(u',v_1)}}{2}\\
&=\frac{2^{d(u',v_1)}}{4}>0.
\end{align*}

Which is a contradiction to $s_\omega(v_0)=s_\omega(T)$, and we are done.
\end{pf}

\begin{cor}
Let $\omega$ be a nonnegative function in $V(T)$, for some $v\in W$, $\omega'$ be a nonnegative function satisfies $\omega'(v)=\omega(v)-1$, $\omega'(u)=\omega(u)$ for other vertices in $T$, then
$$s_\omega(T)\geq s_{\omega'}(T)+2^{d_\omega(v,l)}.$$
\label{cor1}\end{cor}

\begin{pf}
Assume that there exist $v_1$ and $v_2$, so that $s_\omega(v_1)=s_\omega(T)$ and $s_{\omega'}(v_2)=s_{\omega'}(T)$.

By the definition of $s_\omega(v)$, if $\omega(v)\geq2$, then $d_\omega(v,l)=d_{\omega'}(v,l)$, we have
\begin{align*}
s_\omega(T)=s_{\omega}(v_1)&\geq s_{\omega}(v_2)\\
&=s_{\omega'}(v_2)+2^{d(v,v_2)}\\
&\geq s_{\omega'}(v_2)+2^{d_{\omega'}(v,l)}~~{\rm (by~Lemma~\ref{lem5})}\\
&=s_{\omega'}(T)+2^{d_\omega(v,l)}.
\end{align*}

If $\omega(v)=1$, the difference between $\overrightarrow{T}\backslash T_\omega(v_1)$ and $\overrightarrow{T}\backslash T_{\omega'}(v_2)$ is just the length of the maximal path containing $v$, so we have
\begin{align*}
s_\omega(T)=s_{\omega}(v_1)&\geq s_{\omega}(v_2)\\
&=s_{\omega'}(v_2)+2^{d(v,v_2)}+2^{d_\omega(v,l)}-2^{d_{\omega'}(v,l)}\\
&\geq s_{\omega'}(v_2)+2^{d_\omega(v,l)}~~{\rm (by~Lemma~\ref{lem5})}\\
&=s_{\omega'}(T)+2^{d_\omega(v,l)}.
\end{align*}
\end{pf}

\begin{thm}
Let $\omega$ be a nonnegative function on $V(T)$, the $\omega$-cover pebbling number of $T$ is
$$\gamma_\omega(T)=s_\omega(T).$$
\label{thm2}\end{thm}

\begin{pf}
The lower bound holds clearly, for we put $2^{a_i}-1$ pebbles on the leaf of each path for $1\leq i\leq n$ (then no pebble can be moved to $T_\omega(v)$), and $\sum_{u\in S}w(u)2^{d(u,v)}-1$ pebbles on $v$, obviously it is not $\omega$-solvable, and we are done.

For the upper bound, it holds if $|\omega|=1$ or $|W|=1$ by the proof of Theorem~\ref{thm0}, also it holds for $|T|\leq2$ by Theorem~\ref{t1} and Theorem~\ref{t2}. So we may assume that $|\omega|\geq2$, $|W|\geq2$ and $|T|\geq3$.

If the result is false for some $T$ and $\omega$, then we choose one counterexample $T$ and its weight $\omega$, so that $|T|$ and $|\omega|$ are both minimal, that means the upper bound holds for $T'$ and its weight $\omega'$ if $|T'|<|T|$ or $|\omega'|<|\omega|$.

Let $D$ be a distribution on $T$ which is not $\omega$-solvable with $s_\omega(T)$ pebbles, by Theorem~\ref{thm1}, we may assume that all pebbles are distributed on the leaves of $T$.

Let $s_\omega(v_0)=s_\omega(T)$ , there exist $x\in W\backslash v_0$ satisfies $d_{T_\omega(v_0)}(x)=1$, then if $d_T(x)\neq 1$, we can get $d(x,l)>0$, and there exist a nonempty connected component in $T\backslash E(T_\omega(v_0))$ which is connected with $x$, say $T_1$, and $b_1\geq b_2\geq \ldots \geq b_m$ is the size of the maximum path partition of $T_1$.

Case $1$. $D(T_1)$ cannot move a pebble to $x$, then $|D(T_1)|\leq \sum_{i=1}^{m}2^{b_i}-m$,  then we consider $D$ on $T\backslash T_1$, $|D(T\backslash T_1)|\geq s_\omega(T)-D(T_1)\geq s_\omega(T\backslash T_1)$, and also $D(T\backslash T_1)$ is not $\omega$-solvable, a contradiction to the minimum of $|T|$.

Case $2$. $D(T_1)$ can move one pebble to $x$, then it cost us at most $2^{b_1}=2^{d_\omega(x,l)}$ pebbles on $T_1$, the left pebbles on $T$ is not $\omega'$-solvable ($\omega'$ satisfies $\omega'(x)=\omega(x)-1$, and unchanged for other vertices in $T$), so from the minimum of $|\omega|$ and Corollary~\ref{cor1}, we have $|D|< s_{\omega'}(T)+2^{d_\omega(x,l)}\leq s_\omega(T)$, a contradiction to $|D|=s_\omega(T)$.

So we may assume $d_T(x)=1$.

We claim that $D(x)=0$. Otherwise,  let $\omega'$ satisfies $\omega'(x)=\omega(x)-1$ and $\omega'(v)=\omega(v)$ for $v\neq x$. Ignore one pebble already on $x$, we know that $|D|-1$ other pebbles cannot solve $\omega'$, from the minimum of $|\omega|$, we have $|D|-1\leq s_{\omega'}(T)-1$. By Corollary~\ref{cor1}, $s_{\omega'}(T)+1\leq s_\omega(T)$, so $|D|\leq s_\omega(T)-1$, a contradiction to $|D|=s_\omega(T)$, so $D(x)=0$.

%

Assume that $x'\sim x$ in $T$, then we delete $x$, let $C'(x')=C(x')+2C(x)$, and $C'(v)=C(v)$ otherwise. Note that all pebbles are distributed on the leaves of $T$, so $C'(x)=D(x')-\omega(x')-2(D(x)-\omega(x))=-\omega(x')-2\omega(x)$. By Lemma~\ref{lem2}, $D$ is not $\omega$-solvable in $T$ is equivalent to $D$ is not $\omega'$-solvable in $T\backslash x$, where $\omega'(x')=\omega(x')+2\omega(x)$ and $\omega'(v)=\omega(v)$ for $v\neq x$. By the minimum of $|T|$, we have $|D|\leq s_{\omega'}(T\backslash x)-1$, note that $x\neq v_0$, we have $s_{\omega'}(T\backslash x)=s_\omega(T)$, a contradiction to $|D|=s_\omega(T)$, and we are done.
\end{pf}


Moreover, from Theorem~\ref{thm2}, we can immediately get
\begin{cor}
Let $T$ be a tree, $\omega$ be a nonnegative function on $V(T)$, $W=\{v\in V(T): \omega(v)>0\}$, $L=\{v\in V(T): d(v)=1\}$, then if $L\subseteq W$,
$$\gamma_\omega(T)=\max_{v\in V(T)}\sum_{u\in V(T)}\omega(u)2^{d(u,v)}.$$
\label{cor2}\end{cor}

Theorem~\ref{t2} gives a sufficient condition of a nonnegative weight function $\omega$ on $V(G)$ for a graph $G$ so that the $\omega$-cover pebbling number of $G$ is
$$\gamma_\omega(G)=\max_{v\in V(G)}\sum_{u\in V(G)}\omega(u)2^{d(u,v)}.$$

Corollary~\ref{cor2} gives a weaker sufficient condition of a nonnegative weight function $\omega$ on $V(T)$ for a tree $T$ so that the $\omega$-cover pebbling number of $T$ is $$\gamma_w(T)=\max_{v\in V(T)}\sum_{u\in V(T)}\omega(u)2^{d(u,v)}.$$
 Here we give some problems.

\begin{pro}
Give a weaker sufficient condition of a nonnegative function $\omega$ on $V(G)$ for a graph $G$ so that the $\omega$-cover pebbling number of $G$ is
$$\gamma_\omega(G)=\max_{v\in V(G)}\sum_{u\in V(G)}\omega(u)2^{d(u,v)}.$$
\end{pro}

\begin{pro}
For a nonnegative function $\omega$, determine the $\omega$-cover pebbling number of more graphs, such as cycles, hypercubes and so on.
\end{pro}

Also we give a conjecture which is similar to Graham's Conjecture.
\begin{cjt}
Let $\omega_1$ be a nonnegative function on $G$, and $\omega_2$ be a nonnegative function on $H$, the function $\omega$ on $G\times H$ is given by $\omega((g,h))=\omega_1(g)\omega_2(h)$, where $g\in V(G)$ and $h\in V(H)$, then $\gamma_\omega(G\times H)\leq \gamma_{\omega_1}(G)\gamma_{\omega_2}(H)$.
\end{cjt}
%
%
%

\textbf{Acknowledgments.} The authors are grateful for the many useful comments provided referees.

\bibliographystyle{plain}
\bibliography{cpn20190628}
\end{document}